\documentclass[12pt]{article}
\usepackage{amsmath,amssymb,amsthm}
\usepackage{graphicx}
\newenvironment{Proof}{{\it Proof.}}{{$\square$} \vskip.5cm}
\newtheorem{Thm}{Theorem}[section]
\newtheorem{Lem}[Thm]{Lemma}

\theoremstyle{definition}

\author{P. Svetlov\footnote{Partially supported by RFBR (grants nos.
02-01-00090) and by the grant NSh-1914.2003.1}}
\title{Obstructions for generalized graphmanifolds\\ to be nonpositively curved}
\date{December 13, 2004}
 \theoremstyle{remark}

\newcommand{\R}{\mathbb{R}}
\newcommand{\Z}{\mathbb{Z}}

\newcommand{\0}{{\bf 0}}

\newcommand{\vp}{\varphi}

 \makeatletter
\renewcommand{\section}{\@startsection{section}{1}{\parindent} %
{5.5ex plus 1ex minus .2ex}{1ex plus .2ex}{\bf\large\!\!\!\!\!\!}}
\renewcommand{\subsection}{\@startsection{subsection}{2}{\parindent} %
{.5ex plus 1ex minus .2ex}{0ex plus .2ex}{\bf\normalsize}}
\renewcommand{\@eqnnum}{{\rm ($\Phi$\theequation\rm)}}              
\renewcommand{\p@equation}{$\Phi$}                                  
 \makeatother

\begin{document}
\maketitle

\begin{abstract}
An $n$-dimensional manifold $M$ ($n\ge 3$) is called {\it
generalized graph manifold}  if it is glued of blocks that are
trivial bundles of $(n-2)$-tori over compact surfaces (of negative
Euler characteristic) with boundary.  In this paper two
obstructions for generalized graph manifold to be nonpositively
curved are described.

Each 3-dimensional generalized graph manifold with boundary
carries a metric of nonpositive sectional curvature in which the
boundary is flat and geodesic (B. Leeb). The last part of this
paper contains an example of 4-dimensional generalized graph
manifold with boundary, which does not admit a metric of
nonpositive sectional curvature with flat and geodesic boundary.
\end{abstract}

\section{Introduction}

\subsection{Generalized graph manifolds.}~ Let $n\ge 3$. A
generalized graph manifold is a compact $n$-di\-men\-si\-o\-nal
manifold consisting of finitely many blocks
$$
M=\bigcup_{v\in V}M_v.
$$
We require that the following conditions 1-3 to be fulfilled.
\begin{enumerate}
\item Each block $M_v$, $v\in V$  is homeomorphic to
$T^{n-2}\times S_v$, where $T^{n-2}$ is the $(n-2)$-dimensional
torus and $S_v$ is a compact oriented surface with boundary and
with negative Euler characteristic. \item The manifold $M$ is
glued from its blocks by diffeomorphisms between their boundary
$(n-1)$-tori. \item The gluing diffeomorphisms do not identify the
homotopy classes of the fiber $(n-2)$-tori.
\end{enumerate}
Note that the decomposition is unique (up to isotopy). If $n=3$
our definition gives a class $\mathfrak M_0$ of irreducible Haken
graph manifolds. Each irreducible Haken graph manifold $M$ has a
covering manifold $\widetilde{M}\in \mathfrak M_0$ (see [1]).

\subsection{Geometrization of a block.}~Let $M_v\simeq T^{n-2}\times S_v$ be a block in a
$n$-dimensional generalized graph manifold $M$. A metric of
nonpositive sectional curvature on $M_v$ is said to be {\it
geometrization} of $M_v$ if
\begin{enumerate}
\item all boundary $(n-1)$-tori of $M_v$ are flat and geodesic;
\item all fiber $(n-2)$-tori in $M_v$ are flat, geodesic and
isometric.
\end{enumerate}
A block  is geometrizable if and only if it admits a metric of
nonpositive sectional curvature in which its boundary tori are
flat and geodesic (see [2], and  [3]).

\subsection{Geometrization of a manifold.}~ We say that a generalized graph manifold $M$ is {\it
geometrizable} if it admits a metric of nonpositive sectional
curvature which induced a geometrization on each its block. A
generalized graph manifold $M$ is geometrizable if and only if and
only if it admits a metric of nonpositive sectional curvature in
which its boundary tori are flat and geodesic (if $M$ is closed,
then it is geometrizable iff it admits a metric of nonpositive
sectional curvature, see [2], and [3]).

In [2], S. V. Buyalo and V. L. Kobel'skii found  necessary and
sufficient conditions for closed generalized graph manifolds to
carry a metric of nonpositive sectional curvature.

\subsection{Results.}~ In this paper two
topological invariants $c\,(M), c'\,(M)\in\Z$ of generalized graph
manifold $M$ are described (see \ref{ct} and \ref{c't}). If any of
the numbers is nonpositive then $M$ does not admit a
geometrization.

Each 3-dimensional generalized graph manifold with boundary admits
a geometrization (B. Leeb). The last section contains an example
of 4-dimensional generalized graph manifold with boundary, which
does not admit a geometrization (theorem \ref{example}).

\paragraph{Acknowledgements.} The author is grateful to professor
S. Buyalo for useful discussions.

\section{Some notations} Let $M$ be a generalized gaph manifold
and let $V$ be the set of its blocks. Let $v\in V$. We define a
set $\partial v$ as the set of boundary components of the block
$M_v$:
$$
\partial M_v=\bigcup_{w\in\partial v}T_w
$$ (each $T_w$ is $(n-1)$-dimensional torus).
Now we divide all boundary tori in two classes:
$$W=\left\{w\in\bigcup_{v\in V}\partial v\,|\,T_w\not\subset\partial
M\right\}, \quad W_\partial=\left\{w\in\bigcup_{v\in V}\partial
v\,|\,T_w\subset\partial M\right\}\,.
$$

A torus $T_w$, $w\in W$ is called {\it gluing torus}. The gluing
diffeomorphisms induce a permutation $\vp: W\to W$ (the torus
$T_w$, $w\in W$ is glued with torus $T_{\vp(w)}$ in $M$). If $w\in
W$, we use the notation $-w$ for $\vp(w)$ and $\psi_w:T_w\to
T_{-w}$ for the gluing diffeomorphism.

\section{Geometrization of a block}

\subsection{Waldhausen bases.}~ Let $L_w=H_1(T_w,\Z)\simeq\Z^{n-1}$, $w\in\partial v$
and $\iota_v^w: F_v\to L_w$ be the inclusion of the first homology
group of a fiber of the block $M_v=T^{n-2}\times S_v$. Note, that
$\iota_v^w(F_v)=F_w$ is a maximal subgroup of $L_w$.  A set
$\{z_w\in L_w|\,w\in
\partial v\}$  {\it represents a boundary} if
$$
\sum_{w\in \partial v}(\iota_w)_*z_w=0\in H_1(M_v,\Z),
$$
where $\iota_w:T_w\to M_v$ is the inclusion 
$T_w\subset \partial M_v\subset M_v$.

Let $(f_v^1,\ldots,f_v^{n-2})$ be a basis in $F_v$ and let a set
$\{z_w\in L_w|\,w\in \partial v\}$ { represents a boundary}. Then
for each $w\in\partial v$ the set
$$(z_w,f_w)=(z_w,\iota^w_vf_v^1,\ldots,\iota^w_vf_v^{n-2})$$ is a
basis in $L_w$. 
The set $\{(z_w,f_w)\,|\,w\in\partial v\}$ is called {Waldhausen
basis} of the block $M_v$. If
$(\bar{f}_v^1,\ldots,\bar{f}_v^{n-2})$ is another basis in $F_v$
and a set $\{\bar{z}_w\in L_w|\,w\in \partial v\}$ { represents a
boundary}, then there exists  a  set of matrices $$h_w=\left(\begin{array}{cc} \det{\sigma_v}&0\\
n_w& \sigma_v\end{array}\right)\in SL_{n-1}(\Z),\quad w\in\partial
v$$ such that
$$
(\bar{z}_w,\bar{f}_w)=(z_w,f_w)\cdot h_w = (\det{\sigma_v}\cdot
z_w+f_wn_w,f_w\sigma_v)
$$
for each $w\in\partial v$ (here
$\bar{f}_w=(\iota_v^w\bar{f}_v^1,\ldots,\iota_v^w\bar{f}_v^{n-2})$).
It is clear that $\sum_{w\in\partial v}n_w=0$ (see [2]).

\subsection{The Leeb conditions.}\label{gb}~ Each geometrization of $M_v$ defines (and is defined by)
a set of positively definite quadratic forms $\{\beta_w:
L_w\to\R\,|\,w\in\partial v\}$, such that
\begin{itemize}
\item[L1)] the map $\beta_w\circ\iota_v^w: F_v\to\R$ does not
depend on $w\in\partial v$;

\item[L2)] if a set of classes $\{z_w\in L_w|\,w\in \partial v\}$
{ represents a boundary} then $$\sum_{w\in \partial
v}B_w(z_w,\iota^w_vh)=0\quad \mbox{for any} \quad h\in F_v$$
\end{itemize}
(here $B_w$ is the bilinear form associated with $\beta_w$). The
 conditions L1,2 are called {\it the Leeb conditions} (see [2], and
[3]). We note that  the condition L2 is correct: if it is true for
some set $\{z_w\in L_w\}$ represented a boundary then so is for
each set $\{\bar{z}_w\in L_w\}$ represented a boundary:
$$
\sum_{w\in \partial v}B_w(\bar{z}_w,\iota^w_vh)=\sum_{w\in
\partial v}B_w( \det{\sigma_v}\cdot z_w+f_wn_w,\iota^w_vh)=$$$$=
(\det{\sigma_v})\sum_{w\in
\partial v}B_w({z}_w,\iota^w_*h)+B_v(f_v,h)\sum_{w\in
\partial v}n_w=0.
$$
Here $B_v=(\iota_v^w)^*B_w$ (by L1 it is correctly defined).

\subsection{The Gram matrix.}~ Let
$(z,f)_v=\{(z_w,f_w)\,|\,w\in\partial v\}$ be a Waldhausen basis
of the block $M_v$. We define a set of {\it Gram matrices}
$\{G_w^{(z,f)_v}\,|\,w\in\partial v\}$ for the geometrization
$\{\,\beta_w\}_{w\in\partial v}$ in the basis $(z,f)_v$ as
follows:
$$
G_w^{(z,f)_v}=\left(\begin{array}{ccc} x_w& &l_w\\
&& \\l_w^T& & G_v\end{array}\right),
$$
where the number $x_w$, the $(n-2)$-row $l_w$, and the square
matrix $G_v$ are defined by $$x_w=B_w(z_w,z_w),\quad
(l_w)_i=B_w(z_w,\iota_v^wf_v^i),\quad
(G_v)_{ij}=B_v(f_v^i,f_v^j)\,.$$ The condition L2 gives
$$\sum_{w\in\partial v}l_w=\sum_{w\in\partial v}B_w(z_w,\iota_v^wf_v)=0,$$
so the row $\sum_wl_w$ is zero. It is clear, that if
$(\bar{z}_w,\bar{f}_w)=(z_w,f_w)\cdot h_w$ is another Waldhausen
basis, then $ G_w^{(\bar{z},\bar{f})_v}=h_w^TG_w^{(z,f)_v}h_w. $

\section{The first obstruction}

Let $M$ be a n-dimensional generalized graph manifold, let $M_v$,
$M_{v'}$ be blocks in $M$ and  let $\psi_w:T_w\to T_{-w}$ be the
gluing map (here $w\in\partial v$, $-w\in\partial v'$). If the
generalized graph manifold $M$ has a geometrization (as it
described in \ref{gb}) then the induced maps $\Psi_w:
(L_w,\beta_w)\to (L_{-w},\beta_{-w})$
$$
(z_{-w},f_{-w})=(z_w,f_w)\Psi_w
$$
are isometries.
 Indeed, the homogeneous system  of linear
equations
$$
\left\{\begin{array}{cc} G_{-w}=\Psi_w^TG_w\Psi_w,&w\in W\\
\sum_{w\in \partial v}l_w=0,&v\in
V\end{array}\right.\eqno{LS_1(M)}
$$ has a nontrivial solution
$$
G_w=\left(\begin{array}{ccc} x_w& &l_w\\
&& \\l_w^T& & G_v\end{array}\right),\quad w\in \partial v,\quad
v\in V.
$$
Let $c\,(M,(z,f))$ --- be corank of the linear system (the number
of its variables minus its rank). Note that a homogeneous  linear
system has a nontrivial solution iff its corank is positive.

\begin{Thm}\label{ct}
The integer $c\,(M,(z,f))$ does not depend on the choice of the
basis $(z,f)$.
\end{Thm}
\begin{Proof} If we change Waldhausen bases in blocks of $M$ then
the matrix $A$ of $LS_1(M)$ changes to a matrix $AH$, where $H$ is
nonsingular matrix.\end{Proof}

So $c\,(M)=c\,(M,(z,f))$ is a topological invariant of $M$.
 Note, that if $c(M)\le 0$ then $M$ does not admit
geometrization.
\begin{Lem}\label{cl}
If $M$ is a $n$-dimensional generalized graph manifold $M$ then
   $$
   c(M)\ge
   (n-1)|W_\partial|+\frac{(n-2)(n-3)}{2}\cdot|V|-\frac{(n-1)(n-4)}{4}\cdot|W|
   .
   $$
\end{Lem}\begin{Proof} By direct calculation with using the obvious
inequality between the rank and the number of equations in the
linear system.\end{Proof}

If $n=3,4$ then $c(M)>0$ for any $n$-dimensional generalized graph
manifold $M$. So the obstruction to geometrization vanishes.

\section{The second obstruction}
 Let $M$ be a n-dimensional generalized graph manifold,
let $M_v$, $M_{v'}$ be blocks in $M$ and  let $\psi_w:T_w\to
T_{-w}$ be the gluing map (here $w\in\partial v$, $-w\in\partial
v'$). Assume that $M$ has a geometrization (as it described in
\ref{gb}). Then the induced maps $\Psi_w: (L_w,\beta_w)\to
(L_{-w},\beta_{-w})$
$$
(z_{-w},f_{-w})=(z_w,f_w)\Psi_w
$$
are isometries.

Let $n>3$. Then the intersection
$$I_w=\Psi_{-w}(\iota_{v'}^{-w}F_{v'})\cap \iota_v^{w}F_v$$ is non-empty
$(n-3)$-dimensional sublattice in $L_{w}$. The restriction
$$\Psi_w|_{I_w}:(I_w,\beta_w|_{I_w})\to (I_{-w},\beta_{-w}|_{I_{-w}})$$
must be isometry. Let
$$
\Psi_w=\left(\begin{array}{cc} a_w& b_w\\c_w&
d_w\end{array}\right)
$$
be the matrix in some Waldhausen bases, and let $P_w$ be a nonzero
$(n-2)\times(n-3)$ matrix of rank $n-3$ s.t. $b_wP_w=0$. The
equality $f_{-w}=z_wb_w+f_wd_w$ implies $f_{-w}P_w=f_wd_wP_w$. The
vectors consist a basis of $I_{w}\simeq I_{-w}$, so we have
$$ P_w^T(G_{v'}-d_w^TG_vd_w)P_w=0\eqno{LS_2(M)}
$$
Let $c'(M,(z,f))$ --- be corank of the linear system (the number
of its variables minus its rank). Note that a linear system has a
nontrivial solution iff its corank is positive.

The following two proposition are proved in the same way as
\ref{ct} and \ref{cl}.

\begin{Thm}\label{c't}
The integer $c'(M,(z,f))$ does not depend on the choice of the
basis $(z,f)$.
\end{Thm}

So $c'(M)=c'(M,(z,f))$ is a topological invariant of $M$. Note,
that if $c'(M)\le 0$ then $M$ does not admit geometrization.
\begin{Lem}\label{c'l}
If $M$ is a $n$-dimensional generalized graph manifold $M$ and
$n>3$ then
$$
c'(M)\ge \frac{(n-2)(n-1)}{2}\cdot|V|-\frac{(n-3)(n-2)}{4}\cdot
|W|.
$$
\end{Lem}

\section{The example} Let $M$ be a generalized graph manifold with one boundary torus pasted from one block
 $M_1=S\times T^2$ (here $S$ is
an oriented surface with 7 boundary circles) along gluing maps
$\psi_i: T_i\to T_{-i}$, $i=1,2,3$, where
$$
\partial M_1=\bigcup_{i=-3}^3T_i\,,\quad\mbox{and}\quad
\Psi_i=\left(\begin{array}{ccc} 1& 1& i-1\\0& 1& 1\\0& 3&
2\end{array}\right),\,i=1,2,3
$$
(in some Waldhausen basis).  Note, that $c'(M)\ge 0$ (lemma
\ref{c'l})
\begin{Thm}\label{example}
The generalized graph manifold $M$ with boundary $\partial M=T_0$
does not admit a geometrization.
\end{Thm}
{\it Proof.}  The linear system $LS_2(M)$ is
$$
P_i^T(G_v-d_i^TG_vd_i)P_i=0,\quad i=1,2,3\,\eqno{(*)}
$$
where
$$
d_i=\left(\begin{array}{cc}  1& 1\\ 3& 2\end{array}\right) \,,
\quad\mbox{and}\quad
P_i=\left(\begin{array}{c}i-2\\3i-1\end{array}\right).$$ So, the
system ($*$) can be presented as follows:
$$
(15i^2-52i+45)x+(66i^2-232i+204)y+(72i^2-258i+231)z=0,\,i=1,2,3.
$$
(here $x,y,z$ are entries of $G_v$). The matrix of ($*$) is
$$A=\left(\begin{array}{ccc}
8& 38&45\\1& 4& 3\\24&102& 105\end{array}\right)
$$
and rank of $A$ is 3, so $c'(M)=0$. Therefore the manifold $M$
does not admit a geometrization. $\square$

\begin{itemize}
\item[[1\hspace{-5pt}]] S. Buyalo, V. Kobel'skii, {\it
Geometrization of graphmanifolds. I. Conformal geometrization},
St. Petersburg Math. J., {\bf 7} (1996), No. 2, 185-216

\item[[2\hspace{-5pt}]] S. Buyalo, V. Kobel'skii, {\it Generalized
graphmanifolds of nonpositive curvature},  St. Petersburg Math.
J., {\bf 11} (2000), No. 2, 251-268

\item[[3\hspace{-5pt}]] B. Leeb, {\it 3-manifolds with(out)
metrics of nonpositive curvature}, Invent. Math., {\bf 122},
(1995), 277-289
\end{itemize}

\vspace{1cm}

{\tt St. Petersburg Branch,

Steklov Mathematical Institute,

Fontanka 27, 191023 St.Petersburg, Russia}

\vspace{0.5cm} svetlov@pdmi.ras.ru

\end{document}